\documentclass{article}
\usepackage{amsmath,amsthm}
\usepackage{amssymb,yhmath}
\usepackage{yfonts,mathrsfs}
\usepackage{color}

\usepackage{latexsym}
\usepackage[english]{babel}
\usepackage[dvips]{graphicx}

\begin{document}

\newcommand{\Lam}{\Lambda}
\newcommand{\Ga}{\Gamma}
\newcommand{\lam}{\lambda}
\newcommand{\al}{\alpha}
\newcommand{\eps}{\varepsilon}
\newcommand{\de}{\delta}
\newcommand{\om}{\omega}
\renewcommand{\th}{\theta}
\newcommand{\ga}{\gamma}

\newcommand{\R}{\mathbb{R}}
\newcommand{\C}{\mathbb{C}}
\newcommand{\Z}{{\mathbb Z}}
\newcommand{\Zero}{{\mathcal{Z}}}

\newcommand{\bi}{\begin{itemize}}
\newcommand{\ei}{\end{itemize}}
\newcommand{\bd}{\begin{description}}
\newcommand{\ed}{\end{description}}
\newcommand{\be}{\begin{enumerate}}
\newcommand{\ee}{\end{enumerate}}

\renewcommand{\i}{\item}

\newcommand{\bqn}{\begin{eqnarray}}
\newcommand{\eqn}{\end{eqnarray}}
\newcommand{\eqnn}{\nonumber\end{eqnarray}}
\newcommand{\eqnl}[1]{\label{#1}\end{eqnarray}}
\newcommand{\bcas}{\begin{cases}}
\newcommand{\ecas}{\end{cases}}

\newcommand{\nn}{\nonumber\\}
\newcommand{\ba}[1]{\begin{array}{#1}}
\newcommand{\ea}{\end{array}}

\newcommand{\funz}[5]{#1 : \begin{tabular}{ccl}
 #2 &$\rightarrow$& #3, \\
 #4 & $\mapsto$& #5   \end{tabular}}

\newcommand{\mfunz}[5]{\funz{$#1$}{$#2$}{$#3$}{$#4$}{$#5$}}
\newcommand{\mmfunz}[5]{\begin{center}
\funz{$\displaystyle{#1}$}{$\displaystyle{#2}$}{$\displaystyle{#3}$}{$\displaystyle{#4}$}{$\displaystyle{#5}$}
\end{center}}

\newcommand{\fz}[3]{#1:\, #2 \rightarrow #3}

\newcommand{\immagine}[3][8]{
\begin{figure}[hbpt]
\begin{center}
\includegraphics[width=#1cm]{immagini/#2}
\caption{#3}
\label{fig:#2}
\end{center}
\end{figure}}

\newtheorem{ml}{\bf Lemma}
\newtheorem{Theorem}{\bf Theorem}
\newtheorem{mo}{\bf \underline{{\sl Observation}}}
\newtheorem{mrem}{\bf \underline{{\sl Remark}}}
\newtheorem{mcc}{\bf Corollary}
\newtheorem{Definition}{\bf Definition}
\newtheorem{mpr}{\bf Proposition}
\newtheorem{mproperty}{\bf Property}

\newcommand{\bt}{\begin{Theorem}}
\newcommand{\et}{\end{Theorem}}
\newcommand{\bl}{\begin{ml}}
\newcommand{\el}{\end{ml}}
\newcommand{\bp}{\begin{mpr}}
\newcommand{\ep}{\end{mpr}}
\newcommand{\bc}{\begin{mcc}}
\newcommand{\ec}{\end{mcc}}
\newcommand{\bdeff}{\begin{Definition}}
\newcommand{\edeff}{\end{Definition}}
\newcommand{\brem}{\begin{mrem}\rm}
\newcommand{\erem}{\end{mrem}}
\renewcommand{\proof}{\b{Proof: }}

\newcommand{\Pt}[1]{\left( #1 \right)}
\newcommand{\Pg}[1]{\left\{ #1 \right\}}
\newcommand{\Pq}[1]{\left[ #1 \right] }
\newcommand{\Pa}[1]{\left\langle #1 \right\rangle }
\newcommand{\virg}[1]{\textquotedblleft#1\textquotedblright}

\renewcommand{\r}[1]{(\ref{#1})}
\newcommand{\con}{\mathit{C}}
\renewcommand{\l}[1]{\mathscr{L}\Pt{#1}}
\renewcommand{\b}[1]{{\bf{#1}}}
\newcommand{\ffoot}[1]{\footnote{#1}}
\newcommand{\dom}{\mathscr{D}}
\newcommand{\dir}[1]{\mathrm{Dir}\Pt{#1}}
\newcommand{\Cs}{\mathscr{C}}
\renewcommand{\H}{\mathscr{H}}
\newcommand{\Id}{\mathrm{Id}}
\newcommand{\sign}{\mathrm{sign}}

\newcommand{\p}{\b{p}}
\newcommand{\bth}{\b{\th}}
\newcommand{\g}{\b{g}}
\newcommand{\nuovo}[1]{{#1}}


\thispagestyle{empty}
\begin{center} \noindent
{\LARGE{\sl{\bf \nuovo{Existence of planar curves minimizing length and curvature}}}}
\vskip 1cm
Ugo Boscain

{\footnotesize CNRS CMAP, \'{E}cole Polytechnique,  Route de Saclay,  91128 Palaiseau Cedex, France}

and

{\footnotesize SISSA, via Beirut 2-4 34014 Trieste, Italy - {\tt boscain@sissa.it}}\\
\vspace*{.2cm}

Gr\'{e}goire Charlot\\
{\footnotesize Institut Fourier - UMR5582, Institut Fourier, 100 rue des Maths, BP 74, 38402 St Martin d'Heres, France - {\tt Gregoire.Charlot@ujf-grenoble.fr}}\\
\vspace*{.2cm}

Francesco Rossi\\
{\footnotesize SISSA, via Beirut 2-4 34014 Trieste, Italy - {\tt rossifr@sissa.it}}

\end{center}
\vspace{.2cm} \noindent \rm
\begin{quotation}
\noindent  {\bf Abstract}
In this paper we consider the problem of reconstructing a curve that is partially hidden or corrupted by minimizing the functional $\int \sqrt{1+K_\ga^2} \,ds$, depending both on length and curvature $K$. We fix starting and ending points as well as initial and final directions.

For this functional we discuss the problem of existence of minimizers on various functional spaces. We find non-existence of minimizers in cases in which initial and final directions are considered with orientation. In this case, minimizing sequences of trajectories can converge to curves with angles.

We instead prove existence of minimizers for the ``time-reparameterized'' functional $$\int \| \dot\gamma(t)  \|\sqrt{1+K_\ga^2}
 \,dt$$  for all boundary conditions if initial and final directions are
considered regardless to orientation. In this case, minimizers can
present cusps (at most two) but not angles.
\end{quotation}



\section{Problems statements and main results}
\label{s:planprobs}

Consider a smooth function $\fz{\ga_0}{[a,b]\cup[c,d]}{\R^2}$ (with $a<b<c<d$) representing a curve that is partially hidden or deleted in $(b,c)$. We want to find a curve $\fz{\ga}{\Pq{b,c}}{\R^2}$ that completes $\ga_0$ in the deleted part and that minimizes a cost depending both on length $\l{\ga}$ and curvature $K_\ga$.

The fact that $\ga$ completes $\ga_0$ means that $\ga(b)=\ga_0(b),\ \ga(c)=\ga_0(c)$. It is also reasonable to require that the directions of tangent vectors (with orientation) coincide, i.e. $\dot{\ga}(b)\sim\dot{\ga_0}(b),\ \dot{\ga}(c)\sim\dot{\ga_0}(c)$ where 
\bqn
v_1\sim v_2 \mbox{ if it exists }\al\in\R^+\mbox{ such that } v_1=\al\, v_2.
\eqnl{intro:ident-verso}
We call these conditions \b{boundary conditions with orientation}.
All along this paper we assume that starting and ending points never coincide, i.e. $\ga_0(b)\neq\ga_0(c)$, and that initial and final directions are nonvanishing.

In the literature this problem has been deeply studied for its application to problems of segmentation of images (see e.g. \cite{bellettini,masnou,linner1,linner2}) and for the construction of spiral splines \cite{coope}. 

The cost studied in \cite{masnou,coope,linner1} is the total squared curvature $E_1\Pq{\ga}=\int_0^{\l{\ga}} |K_\ga(s)|^2\, ds$  where $s$ is the arclength. In \cite{coope,linner1} boundary conditions differ from our boundary conditions with orientations. In particular, starting and ending directions are fixed with angles measured in $\R$ (while we identify $\al$ and $\al +2 k \pi$). In this framework, non existence of minimizers is proved if the starting and ending angles $\th_0,\th_1$ satisfy $|\th_1-\th_0|>\pi$.

The cost studied in \cite{bellettini} is $E_2\Pq{\ga}=\int_0^{\l{\ga}} \Pt{1 + |K_\ga(s)|^2}\, ds$, while in \cite{linner2} it is $E_3\Pq{\ga}=\int_0^{\l{\ga}} \Pt{\eta + |K_\ga(s)|^2}\, ds$ with $\eta\rightarrow 0$. Depending on the cost, minimizers may present angles and the curvature becomes a measure.

\nuovo{The cost $E_4\Pq{\ga}=\int_0^{\l{\ga}} \sqrt{1 + |K_\ga(s)|^2}\, ds$ naturally arises in problems of geometry of vision \cite{citti-sarti,petitot,petitot-libro}.} In this paper we study the following cost:
\bqn
J\Pq{\ga}=\int_b^c \sqrt{\| \dot{\ga}(t)\|^2 + \| \dot{\ga}(t)\|^2 K^2_\ga(t)}\, dt,
\eqnl{intro:costocurva}
that is an extension of $E_4\Pq{\ga}$, see Remarks \ref{rem:Jbello}-\ref{rem:Jbrutto} below. Using this cost one can study the existence of minimizers with angles without involving sophisticated functional spaces. Moreover, this problem has also been studied in \cite{nostro-PTS2}, where it is defined on the sphere $S^2$ instead of the plane.
\brem
The cost $J$ is invariant both by rototranslation and reparametrization of the curve.

Define the cost $J_\beta[\ga]:=\int_b^c \sqrt{\| \dot{\ga}(t)\|^2 + \beta^2 \| \dot{\ga}(t)\|^2 K^2_\ga(t)}\, dt$ with a fixed $\beta\neq 0$. Consider an homothety $(x,y)\mapsto(\beta x,\beta y)$ and the corresponding transformation of a curve $\gamma=(x(t), y(t))$ to $\gamma_\beta=(\beta x(t),\beta y(t))$. It is easy to prove that $J_\beta\Pq{\gamma_\beta}=\beta^2 J\Pq{\gamma}$. Hence the problem of minimization of $J_\beta$ is equivalent to the minimization of $J$ with a suitable change of boundary conditions. Thus, results about $J$ given in this paper hold also for the cost $J_\beta$.
\erem

The first question we address in this paper is the choice of a set of smooth curves on which this cost is well-defined. We want $\dot{\ga}(t)$ and $K_\ga(t)=\frac{\dot x \ddot y - \dot y \ddot x}{\Pt{\dot x^2+\dot y^2}^{\frac{3}{2}}}$ well-defined, thus it is reasonable to look for minimizers in $$\dom_1:=\Pg{\ga\in\con^2([b,c],\R^2)\ \mid\ \dot{\ga}(t)\neq 0\ \forall\, t\in[b,c],\, \ga(b)=\ga_0(b),\, \ga(c)=\ga_0(c),\,
 \dot\ga(b)\sim\dot\ga_0(b),\, \dot\ga(c)\sim\dot\ga_0(c)}.$$
Moreover, we have that $\dot{\ga}(b)$ and $\dot{\ga}(c)$ are well-defined in this case.

\brem
\label{rem:Jbello}
The cost $J\Pq{\ga}$ on the set $\dom_1$ coincide with the cost $E_4\Pq{\ga}$. To prove it, reparametrize a curve in $\dom_1$ by arclength and observe that in this case we have $\|\dot\ga\|=1$.
\erem
Under this assumption, one of the main results of the paper is the nonexistence of minimizers for $J$.
\bp
\label{p:n-ex-1}
There exist boundary conditions $\ga_0(b),\ga_0(c)\in\R^2$ with $\ga_0(b)\neq\ga_0(c)$, $\dot\ga_0(b),\ \dot\ga_0(c)\in\R^2\backslash\Pg{0}$ such that the cost \r{intro:costocurva} does not admit a minimum over the set $\dom_1$.
\ep

To get existence of minimizers for this cost one can choose to enlarge the set of admissible curves. In this paper we consider the simplest generalization, by taking the space 
$$\dom_2:= \Pg{\ga\in\con^2([b,c],\R^2)\ \mid\ \| \dot{\ga}(t)\|^2+ \| \dot{\ga}(t)\|^2 K^2_\ga(t)\in L^1([b,c],\R),\,
\ba{cc}
\ga(b)=\ga_0(b),\,& \ga(c)=\ga_0(c),\\
 \dot\ga(b)\sim\dot\ga_0(b),\,& \dot\ga(c)\sim\dot\ga_0(c)
\ea},$$
on which the cost $J\Pq{\ga}$ is defined and always finite.
\brem
\label{rem:Jbrutto}
Notice that the cost $E_4\Pq{\ga}$ is not well defined on $\dom_2$, since it is not possible in general to perform an arclength parametrization. Then $J\Pq{\ga}$ is an extension of $E\Pq{\ga}$, since they coincide on $\dom_1$.
\erem

Also on $\dom_2$ we have non-existence of minimizers for $J$.
\bp
\label{p:n-ex-2}
There exist boundary conditions $\ga_0(b),\ga_0(c)\in\R^2$ with $\ga_0(b)\neq\ga_0(c)$, $\dot\ga_0(b),\ \dot\ga_0(c)\in\R^2\backslash\Pg{0}$ such that the cost \r{intro:costocurva} does not admit a minimum over the set $\dom_2$.
\ep

The basic problem is that we can have a sequence of minimizing curves converging to a non admissible curve. In particular, we can have angles at the beginning and/or at the end, i.e. each curve $\ga_n$ satisfies given boundary conditions with orientation but the limit curve $\bar\ga$ doesn't satisfy them. See Figure \ref{fig:angoli-inizio}.
\immagine{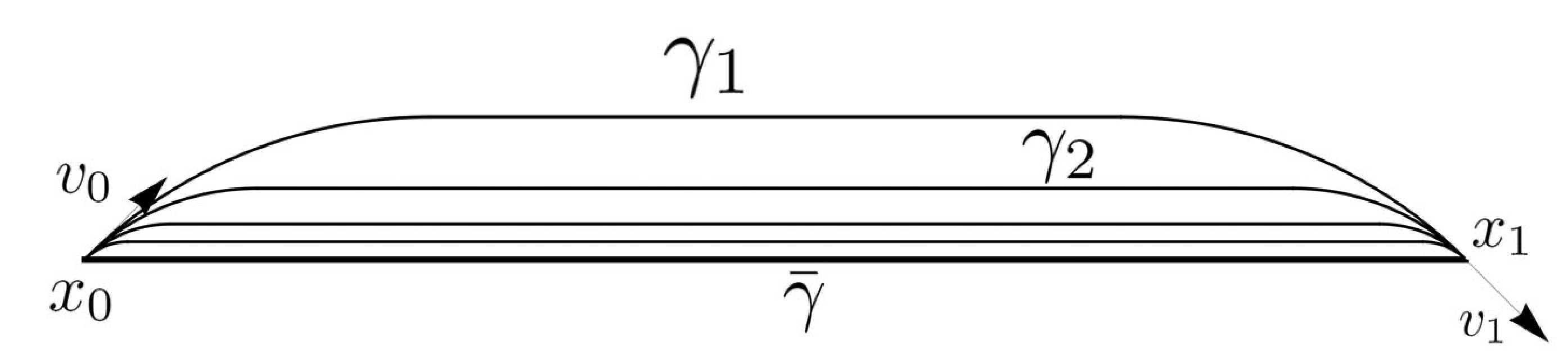}{Minimizing sequence converging to a non-admissible curve (angles at the beginning/end).}



The main result of the paper is the existence of minimizers for the cost \r{intro:costocurva} taking again curves for which $\| \dot{\ga}(t)\|^2+ \| \dot{\ga}(t)\|^2 K^2_\ga(t)$ is integrable, but changing boundary conditions. We only impose conditions on the direction of $\dot{\ga}$ regardless of its orientation.

As before, fix a starting point $x_0$ with a direction $v_0$ and an ending point $x_1$ with a direction $v_1$. Consider planar curves satisfying the following boundary conditions: $\ga(0)=x_0,\ \dot{\ga}(0)\approx v_0,\ \ga(T)=x_1,\ \dot{\ga}(T)\approx v_1$, where the identification rule $\approx$ is 
\bqn
v_1\approx v_2 \mbox{ if it exists }\al\in\R\backslash\Pg{0}\mbox{ such that } v_1=\al v_2.
\eqnl{intro:id-proj}
We call them \b{projective boundary conditions}. As already stated, we have the following existence result:
\bp
\label{p:ex-3}
For all boundary conditions $x_0,x_1\in\R^2$ with $x_0\neq x_1$, $v_0,\ v_1\in\R^2\backslash\Pg{0}$, the cost \r{intro:costocurva} has a minimizer over the set $$\dom_3:=\Pg{\ga\in\con^2([b,c],\R^2)\ \mid\ \| \dot{\ga}(t)\|^2+ \| \dot{\ga}(t)\|^2 K^2_\ga(t)\in L^1([b,c],\R),\, \ga(b)=x_0,\ \ga(c)=x_1,\ \dot\ga(b)\approx v_0,\ \dot\ga(c)\approx v_1}.$$
\ep

Observe that we can have minimizers with cusps, as $\bar\ga$ in Figure \ref{fig:cusp-una}. Indeed, the limit direction (regardless to orientation) is well defined in the cusp point, while the limit direction with orientation is undefined.
\immagine{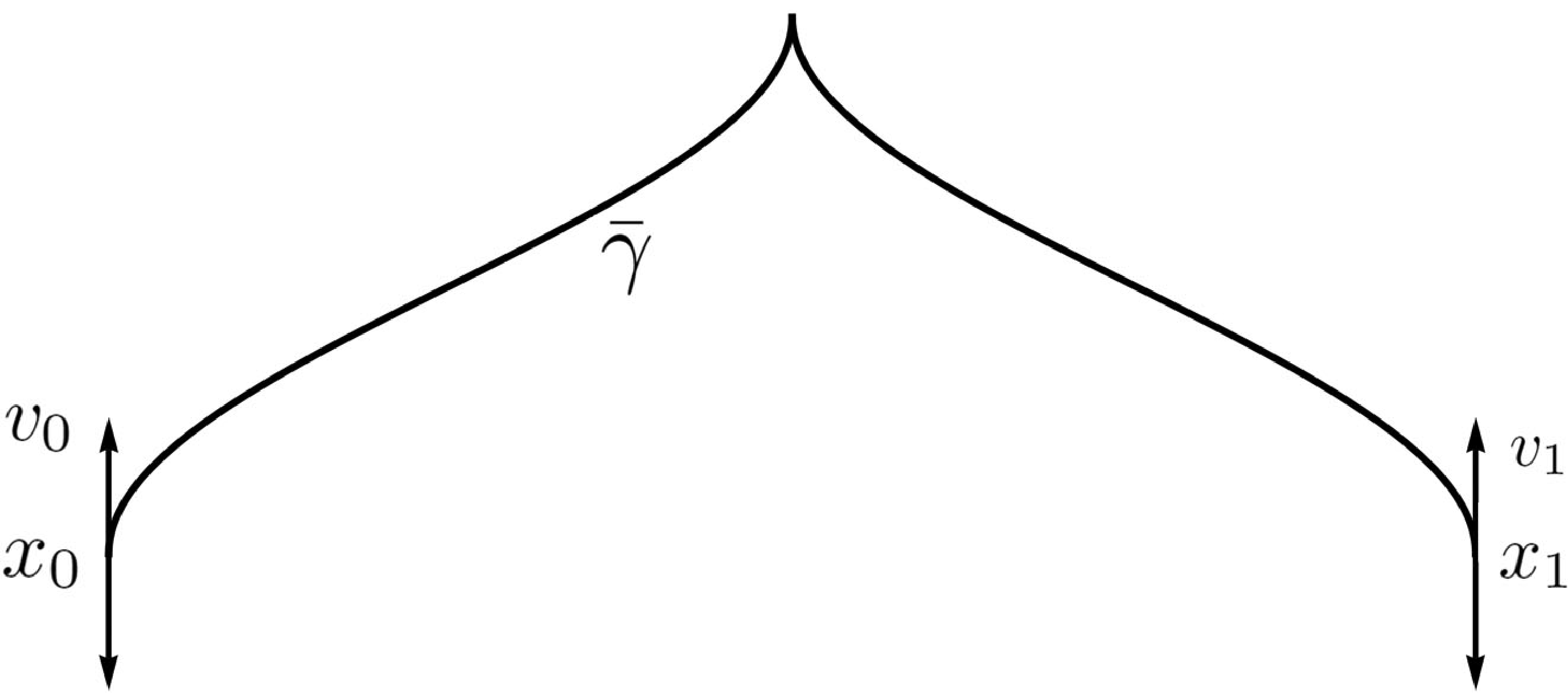}{A minimizer with a cusp.}

All the previous results are obtained as consequences of the study of two similar mechanical problems. For what concerns problems with boundary conditions with orientation, we consider a car on the plane that can move only forwards and rotate on itself (it is the Dubins' car, see \cite{dubins}). Fix two points $(x_0,y_0),\ (x_1,y_1)$ and two angles $\th_0, \th_1$ in these points measured with respect to the positive $x$-semiaxis. Consider all trajectories $q(.)$ steering the car from the point $(x_0,y_0)$ rotated of an angle $\th_0$ to the point $(x_1,y_1)$ rotated of an angle $\th_1$. Our goal is to find the cheapest trajectory with respect to a cost depending both on the length of displacement on the plane and on the angle of rotation on itself.

The dynamics can be written as the following control system on the group of motions of the plane $SE(2):=\Pg{\Pt{x,y,\th}\ |\ (x,y)\in\R^2,\ \th\in\R/2\pi}$:
\bqn
&&\Pt{\ba{c}
\dot x\\
\dot y\\
\dot \th
\ea
}=u_1\Pt{\ba{c}
\cos(\th)\\
\sin(\th)\\
0
\ea}+u_2\Pt{\ba{c}
0\\
0\\
1
\ea}
\eqnl{intro:dynpiano}
where $x,y$ are coordinates on the plane and $\th$ represents the angle of rotation of the car. Since we forbid backwards displacements, we impose $u_1\geq 0$. We want to minimize the cost
\bqn
\Cs\Pq{q(t)}=\int_0^T \sqrt{u_1^2(t)+ u_2^2(t)}\,dt
\eqnl{intro:costocontrollo}
with the following boundary conditions: $x(0)=x_0,\ y(0)=y_0,\ \th(0)=\th_0,\ x(T)=x_1,\ y(T)=y_1,\ \th(T)=\th_1$.

\brem
\label{r:lift} Any smooth planar curve can be naturally transformed into an admissible trajectory of this control system. Indeed, given  $\ga(t)=\Pt{x(t),y(t)},$ we set $q(t)=\Pt{x(.),y(.),\th(.)}$ where $\th(t)$ is the angle of the tangent vector with respect to the positive $x$-semiaxis. This construction is called the \b{lift} of the curve $\ga$. 

We then find suitable controls $u_i$ corresponding to $q(t)$ defined above. In this framework $u_1$ plays the role of $\|\dot{\ga}\|,$ while $u_2$ is $\|\dot{\ga}\|K_\ga$. Hence, the cost \r{intro:costocontrollo} coincides with $J\Pq{\ga}$ defined in \r{intro:costocurva}. Moreover, boundary conditions with orientation can be easily translated to boundary conditions on $(x,y,\th)\in SE(2)$.

Notice that, on the contrary, not all trajectories of \r{intro:dynpiano} are lifts of planar curves. Indeed, consider a trajectory of the system with $u_1\equiv 0$. It represents the rotation of the car on itself. If we consider its projection to the plane $\Pi:(x,y,\th)\mapsto(x,y)$ the curve is reduced to a point, thus $\dot{\ga}=0$ and the curvature is undefined.
\erem

We will prove that for the optimal control problem \r{intro:dynpiano}-\r{intro:costocontrollo} on $SE(2)$ we have existence of minimizers with $L^1$ controls. Starting from a minimizer of this problem, we will find counterexamples to the existence of minimizers of $J$ on $\dom_1$ and $\dom_2$.\\

For what concerns problems with projective boundary conditions, we study the dynamics given by \r{intro:dynpiano} where we admit also backwards displacements (it is the Reeds-Shepp car, see \cite{Reeds-Shepp}). In this case we don't have to impose $u_1\geq 0$ and we identify $(x,y,\th)\simeq(x,y, \th+\pi)$. Hence this dynamics is naturally defined on the quotient space $SE(2)/\simeq$. We choose the same cost \r{intro:costocontrollo}. Also in this case, it is possible to lift planar curves to curves on $SE(2)/\simeq$. Projective boundary conditions can be easily translated to conditions on $(x,y,\th)\in SE(2)/\simeq$.

For the optimal control problem \r{intro:dynpiano}-\r{intro:costocontrollo} on $SE(2)/\simeq$ we have existence of minimizers  with $L^1$ optimal controls. Its consequence on the problem of planar curves is the existence of minimizers of $J$ in $\dom_3$.

For both optimal control problems, the basic tool we use to compute minimizer is the Pontryagin Maximum Principle (PMP in the following), see \cite{pontlibro}. It gives a necessary first-order condition for minimizers. Solutions of PMP are called extremals, hence minimizers have to be found among extremals. For details, see e.g. \cite{agra-book}.\\

The structure of the paper is the following. In Section \ref{s:control} we introduce the group $SE(2)$ and the space $SE(2)/\simeq$, and we define the optimal control problems on these spaces corresponding to the ones defined in Section \ref{s:planprobs} on the plane. We then study the optimal control problems and find some minimizers properties.

Section \ref{s:solution} contains the main results of the paper: we prove Propositions \ref{p:n-ex-1}-\ref{p:n-ex-2}-\ref{p:ex-3} using properties of minimizers of problems studied in Section \ref{s:control}.

\section{Solution of optimal control problems}
\label{s:control}

In this section we recall the definition of the two optimal control problems given above. In the first we consider the Dubins' car \cite{dubins}: it can both move forwards and rotate on itself. In the second we have the Reeds-Shepp car \cite{Reeds-Shepp}, that can move forwards, backwards and rotate on itself. Nevertheless, the problems we study are different than the ones studied in \cite{dubins,Reeds-Shepp}. We don't have constraints on velocity and curvature. We want instead to minimize (in both cases) a cost depending both on velocity and curvature.

\subsection{Dubins' car with lenght-curvature cost}
\label{s:control-or}
The Dubins' car is a car that can move both forwards and rotate on itself. The dynamics of the car is given by the following control system:
\bqn
&&\Pt{\ba{c}
\dot x\\
\dot y\\
\dot \th
\ea
}=u_1\Pt{\ba{c}
\cos(\th)\\
\sin(\th)\\
0
\ea}+u_2\Pt{\ba{c}
0\\
0\\
1
\ea},\qquad u_1\geq 0,
\eqnl{eq:dyncontrollo}
with $u_1,u_2\in L^1([0,T],\R)$. We impose $u_1\geq 0$ to forbid backwards displacements. Observe that $u_1$ is the planar velocity of the car and $u_2$ is its angular velocity. The controllability of this system can be checked by hand, and we omit the proof.

We fix a starting point $q_0=(x_0,y_0,\th_0)$ and an ending point $q_1=(x_1,y_1,\th_1)$. We want to minimize the cost 
\bqn
\Cs\Pq{q(.)}=\int_0^T \sqrt{u_1^2+ u_2^2}
\eqnl{eq:costocontrollo}
over all trajectories of \r{eq:dyncontrollo} steering $q_0$ to $q_1$.  Here the end time $T$ is fixed.
\brem
This problem is a left-invariant problem on the group of motions of the plane $$SE(2):=\Pg{\Pt{
\ba{ccc}
\cos(\th)&-\sin(\th)&x\\
\sin(\th)&\cos(\th)&y\\
0&0&1
\ea}\ |\ (x,y)\in\R^2,\ \th\in\R/2\pi},$$ where the group operation is the standard matrix operation.  Indeed, in this case the dynamics is given by $\dot{g}=u_1\, g p_1 + u_2\, g p_2$ with $$p_1=\Pt{
\ba{ccc}
0&-1&0\\
1&0&0\\
0&0&0
\ea}\qquad p_2=\Pt{
\ba{ccc}
0&0&1\\
0&0&0\\
0&0&0
\ea}.$$ If the constraint $u_1\geq 0$ is removed, we have a minimal length problem on the sub-Riemannian manifold $(SE(2),\Delta,\b{g})$ where $\Delta$ is the left-invariant distribution generated by $p_1,p_2$ at the identity and $\b{g}$ is the metric on $\Delta$ defined in $g$ by the condition $\b{g}_g\Pt{gp_i,gp_j}=\de_{ij}$. For details about sub-Riemannian geometry on Lie groups see e.g. \cite{nostro-gruppi}. For a complete study of this sub-Riemannian problem on $SE(2)$ see \cite{igorino,yuri-solo}.

As a consequence, the problem of minimization of $\Cs$ from $q_0$ to $q_1$ is equivalent to the same problem from $\Id$ to $q_0^{-1} q_1$. For this reason, from now on we will study only problems starting from $\Id$.
\erem

\subsubsection{Existence of minimizers and reduction to $L^\infty$ controls}

\label{s:ex-min-SE2}
In this Section we apply Filippov existence theorem to the optimal control problem \r{eq:dyncontrollo}-\r{eq:costocontrollo} on $SE(2)$, that provides the existence of minimum. We then prove that it is equivalent to solve this problem with controls $u_i\in L^1$ or with controls $u_i\in L^\infty$. This result permits to verify that minimizers found via the PMP, that works in the framework of $L^\infty$ controls, are minimizers also in the larger class of $L^1$ controls.

We first transform problem \r{eq:dyncontrollo}-\r{eq:costocontrollo} in a minimum time problem. It is a standard procedure to transform the problem \r{eq:dyncontrollo}-\r{eq:costocontrollo} with fixed final time $T$ to a problem in which the dynamics is given again by \r{eq:dyncontrollo}, the cost is the time (that is free) and the constraint on the controls are $u_1\geq 0,\ u_1^2+u_2^2\leq 1$.

We apply Filippov existence theorem for minimum time problems, see e.g. \cite[Cor 10.2]{agra-book}, that gives a minimizer, hence $L^1$ optimal controls.

We now prove that we can restrict to $L^\infty$ optimal controls. This generalization cannot be proved in general. Indeed, Lavrentiev phenomenon can occur for more general dynamics and costs, i.e. it may exists a trajectory with $L^1$ controls such that its cost is strictly less that the cost of all trajectories with $L^\infty$ controls, in particular solutions of PMP. For details see e.g. \cite{loewen}.

We thus restrict ourselves to minimal length problems on a trivializable sub-Riemannian manifold with constraints on values of controls.
\bl
\label{l:L1Linf}
Consider a minimal length problem on a trivializable sub-Riemannian manifold $(M,\Delta,\g)$ with constraints on values on control, i.e.
\bqn
\dot q =\sum_{i=1}^m u_i F_i(q)\label{eq:subr}\qquad u(t)\in V\subset \R^n\qquad
\Cs\Pq{q(.)}=\int_0^T\sqrt{\sum_{i=1}^m u_i^2}\ \rightarrow\ \min 
\eqn
where $\Delta(q)=\mathrm{span}\Pg{F_1,\ldots,F_m}$ and $\g_q(F_i(q),F_j(q))=\de_{ij}$. Assume that the set $V$ satisfies $$aV=\Pg{\,a v\ \mid\,\ v\in V}\subset V\mbox{~ ~ ~ for all $a\in \R^+\cup\Pg{0}$}.$$

If there exists a minimizer $\bar{q}(t)$ of this problem with optimal controls $\bar{u}_i\in L^1$, then there exist other optimal controls $\hat{u}_i\in L^\infty$ such that the corresponding trajectory is a minimizer that is a reparametrization of $\bar{q}$.
\el
\proof Let $\bar q:[0,T]\to M$ be a minimizer of the problem \r{eq:subr} with optimal controls in $L^1$. Define 
$$
f(t):=\int_0^t \sqrt{\g\Pt{\dot{\bar q}(\tau),\dot{\bar q}(\tau)}}\,d\tau=\int_0^t\sqrt{\sum_{i=1}^m u_i^2(\tau)}\,d\tau
$$
that is a function from $[0,T]$ to $[0,L]$, with $L=\Cs\Pq{\bar{q}}$. The function $f$ is absolutely continuous and non decreasing, hence the set $R$ of its regular values is of full measure in $[0,L]$. We also define $g:[0,L]\to[0,T]$ by
$$
g(s):=\inf\{t\in[0,T]\;|\;f(t)=s\}.
$$
One can easily check that $\forall\, s\in[0,L]$ it holds $f(g(s))=s$. 

Moreover, if $g$ is discontinuous at $s_0$ then $f^{-1}(s_0)$ is a closed interval of the form $[t_0,t_1]$ and $\forall t\in[t_0,t_1]$, one has $\int_{t_0}^t \sqrt{\g\Pt{\dot{\bar q}(\tau),\dot{\bar q}(\tau)}}\, d\tau=0$, hence $\bar q(t)=\bar q(t_0)=\bar q (g(s_0))$. This also proves that $\bar q\circ g$ is continuous.

We also have that $\bar q\circ g$ is a 1-Lipschtzian function (hence absolutely continuous) since 
$$
d(\bar q(g (s_0)),\bar q(g (s_1)))
\leq \int_{g(s_0)}^{g(s_1)}\sqrt{\g\Pt{\dot{\bar q}(\tau),\dot{\bar q}(\tau)}}\,d\tau
=|s_0-s_1|,
$$
where $d(.,.)$ is the sub-Riemannian distance $d(q_0,q_1):=\inf\Pg{\Cs\Pq{q(.)}\ \mid\ q(.) \mbox{ ~ satisfies \r{eq:subr}~ and steers~ }q_0\mbox{~to~}q_1}$. 

For $s\in R$, $g$ is differentiable at $s$ and its derivative is $\dot g(s)=\frac 1 {\dot f(g(s))}$. We also have that $\bar q$ is differentiable at $g(s)$ because $s$ is a regular value of $f$, which implies that $\dot{\bar q}$ 
is defined. Hence one can easily compute $\g\Pt{\dot{\bar q}(\tau),\dot{\bar q}(\tau)}=1$. Moreover, for $s\in R$ we have
$$
\frac{d(\bar q\circ g)}{ds}(s)=\dot g(s)\dot{\bar q}(g(s))=\frac{1}{\sqrt{\sum_{i=1}^m u_i^2(g(s))}}\dot{\bar q}(g(s)),
$$
hence $\bar q\circ g$ is an admissible curve corresponding to controls 
$$\tilde u_i(s)=\frac{\bar u_i(g(s))}{\sqrt{\sum_{i=1}^m u_i^2(g(s))}},$$ 
that are $L^{\infty}$ controls.

Once these $L^{\infty}$ controls on $[0,L]$ are found, make a linear reparametrization of the time $s\mapsto \frac{s T}{L}$ and a corresponding rescaling of controls $\tilde u_1\mapsto \hat{u}_i:=\frac{\tilde{u}_i L}{T}$. We have now a reparametrization of the same trajectory $\bar{q}$ with controls $\hat{u}_i$ bounded by $\frac{L}{T}$, hence $L^\infty$, on the interval $[0,T]$.
\qed
\brem Our problem \r{eq:dyncontrollo}-\r{eq:costocontrollo} satisfies hypotheses of Lemma \ref{l:L1Linf}, with $V=\Pg{(v_1,v_2)\in\R^2\ |\ v_1\geq 0}$.
\erem
\brem
A consequence of Lemma \ref{l:L1Linf} is that it is equivalent to look for $L^1$ or for $L^\infty$ optimal controls. Indeed, if we have a minimizer with controls in $L^1$, then we reparametrize them and find a minimizer with controls in $L^\infty$. On the opposite side, if we have a minimizer $\bar{q}(t)$ over the set of controls in $L^\infty$, it is also a minimizer over the set of controls in $L^1$. We prove it by contradiction. Let $\tilde{q}(t)$ be a trajectory with controls in $L^1$ and whose cost is less than $\Cs\Pq{\bar{q}(t)}$. Reparametrize $\tilde{q}(t)$ and find a trajectory with controls in $L^\infty$ with the same cost, hence $\bar{q}(t)$ is not a minimizer. Contradiction.
\erem

\subsubsection{Computation of extremals}
\label{s:PMP-SE2}

We now apply the PMP to the problem \r{eq:dyncontrollo}-\r{eq:costocontrollo} transformed into a minimum time problem. For the expression of PMP for minimum time problems see e.g. \cite[Ch. 12]{agra-book}. The control-dependent Hamiltonian of the system is 
\bqn
H(q,\lam,u)=\Pa{\lam,\dot q}=u_1 h_1 + u_2 h_2
\eqnl{eq:Ham}
where $h_1=\lam_x\cos(\th)+\lam_y \sin(\th),\ h_2=\lam_\th$, and $\lam_x,\lam_y,\lam_\th$ are the components of the covector $\lam$ in the dual basis with respect to coordinates $(x,y,\th)$. Notice that $H$ can be seen as the scalar product $(u_1,u_2)\cdot(h_1,h_2)$ in $\R^2$.

We don't give a complete synthesis of the problem, since we only need to find a particular minimizer to use in proofs of Propositions \ref{p:n-ex-1}-\ref{p:n-ex-2}.

We first consider normal extremals, for which we can choose $H=1$. Let us denote $\alpha$ and $\rho$ an angle and a positive number in such a way that $\lambda_x=\rho\cos(\alpha)$ and $\lambda_y=\rho\sin(\alpha)$.

Let us assume that at $t=t_0$ we have $h_1(t_0)>0$. Then PMP gives controls $u_1=h_1$, $u_2=h_2$, after having normalized $\|(h_1,h_2)\|=1$. Dynamics is given by 
\bqn
\begin{cases}
\dot{x}=h_1 \cos(\th)\\
\dot{y}=h_1 \sin(\th)\\
\dot{\lam_x}=\dot{\lam_y}=0\\
\dot{\th}=h_2\\
\dot{\lam_\th}=h_1\Pt{-\lam_x\sin(\th)+\lam_y\cos(\th)}
\end{cases}
\eqnl{eq:dynbella}
We have $|\theta(t_0)-\alpha|<\frac \pi 2$ in $\mathbb{R}/2\pi\mathbb{Z}$ and the equation on $\theta$ is:
$$
2\ddot{\theta}=\rho^2\sin(2(\theta-\alpha)).
$$
It is the equation of the pendulum with $\theta=\alpha$ being the unstable equilibrium. This implies that when starting with $|\theta(t_0)-\alpha|<\frac \pi 2$ in $\mathbb{R}/2\pi\mathbb{Z}$, then $\theta$ will reach a value such that $|\theta(t_0)-\alpha|>\frac \pi 2$ in $\mathbb{R}/2\pi\mathbb{Z}$.
Hence the corresponding extremal will have a time $t_1>t_0$ for which $h_1(t_1)<0$.

Let us assume now that at $t=t_0$ we have $h_1(t_0)\leq 0$. Then PMP gives controls $u_1=0$ and $u_2=\sign(h_2)$. In this case, the extremal corresponds to a rotation on itself. Indeed, the dynamics is given by
\bqn
\begin{cases}
\dot{x}=\dot{y}=\dot{\lam_x}=\dot{\lam_y}=\dot{\lam_\th}=0\\
\dot{\th}=\sign(h_2)\\
\end{cases}
\eqnl{eq:dynbrutta}
Since $\lam_x$ and $\lam_y$ are constant, either they are both vanishing all along the extremal (that is $h_1\equiv 0$) or there exists at least one of them that is nonvanishing. In this case, there exists $t_1>t_0$ such that $h_1(t_1)>0$.
 
As already stated, for an extremal satisfying $h_1(t_0)>0$ (resp $h_1(t_0)<0$) there exists a time $t_1>t_0$ such that $h_1(t_1)<0$ (resp $h_1(t_1)>0$). Thus an extremal is the concatenation of trajectories satisfying \r{eq:dynbrutta}, i.e. pure rotations, and trajectories satisfying \r{eq:dynbella}, that are arcs of pendulum in $\th$.

Consider an arc $\ga(\Pq{t_0,t_1})$ satisfying \r{eq:dynbrutta} between two arcs satisfying \r{eq:dynbella}. Then the variation of $\theta$ along this arc should be of $\pi$ because $\theta(t_0)=\alpha+\frac \pi 2\, \mathrm{mod}\, \pi$ and we should come back to $\theta(t_1)=\alpha+\frac \pi 2\, \mathrm{mod}\, \pi$ at the end with the dynamics $\dot \theta=1$. Moreover, one can prove that a concatenation of dynamics \r{eq:dynbella}, \r{eq:dynbrutta} and \r{eq:dynbella} cannot be optimal.
\brem A consequence of this study to the planar problem of minimization of $J$ on $\dom_2$ is that planar curves with cusps are extremal, but never minimizers. Indeed, a curve with cusp is the projection of a curve $q(.)$ containing an interval in which $\dot x=\dot y=0$, while $\th$ has a variation of $\pi$. Hence, the non optimality of $q$ implies non optimality of the planar curve with cusp.
\erem

\nuovo{
We finally consider abnormal extremals, for which we have $H=0$.
We have two possibilities:
\bi
\i either $h_1=h_2=0$, for which the trajectory is a straight line (i.e. $\dot\th=0$);
\i or $h_2=0, h_1<0$, thus $u_1=0$, for which the trajectory is a pure rotation (i.e. $\dot x=\dot y=0$).
\ei
Abnormal extremals can also be concatenations of these two kind of trajectories.
}

\subsubsection{An example of a minimizer}
\label{s:solstrana}
\newcommand{\q}{\b{q}}
\newcommand{\beps}{{\xi}}

In this section we give an example of a minimizer $\q(.)$ defined on a small interval $\Pq{0,2 \beps}$ and satisfying \r{eq:dynbrutta} on $\Pq{0,\beps}$ and \r{eq:dynbella} on $\Pq{\beps,2\beps}$. This trajectory is the basic example that we will use to prove non-existence of minimizers of cost \r{intro:costocurva} both on $\dom_1$ and $\dom_2$, i.e. Propositions \ref{p:n-ex-1} and \ref{p:n-ex-2}.\\

Consider a trajectory $q^1(t)$ starting from $\Id$, with given $\lam_x=-\frac{1}{\sqrt{2}},\, \lam_y=0,\, \lam_\th=\frac{1}{\sqrt{2}}$. All quantities related to this trajectory are denoted with superscript $1$. Since $h_1(0)<0$, we follow dynamics given by \r{eq:dynbrutta} on an interval $\Pq{0,t^1}$ and we have $$x^1(t^1)=y^1(t^1)=\lam_y^1(t^1)=0,\qquad\lam_x^1(t^1)=-\frac{1}{\sqrt{2}}\qquad \th^1(t^1)=t^1,\qquad \lam_\th^1(t^1)=\frac{1}{\sqrt{2}}.$$ We choose $t^1=\frac{\pi}{2}$ and observe that $h_1^1(t^1)=0$. Recall that on this interval controls are $u_1^1=0,\,u_2^1=1$.

Then dynamics is given by \r{eq:dynbella} on an interval $\Pq{t^1,t^1+s^1}$. Since $\lam_\th^1$ is continuous, so $u_2^1$ is. Then $\th^1(t)=\frac{\pi}{2}+(t-t^1)+o(t-t^1)$ on $\Pt{t^1,t^1+s^1}$, thus $h_1^1(t)>0$ on $\Pt{t^1,t^1+s^1}$ for a sufficiently small choice of $s^1$. As a consequence, $x(t)\neq 0$, $y(t)\neq 0$ for all $t\in\Pt{t^1,t^1+s^1}$, eventually choosing a smaller $s^1$.

Recall now that all normal extremals are local minimizers, i.e. for each extremal $q(t)$  and time $t_0$ there exists $\eps$ such that $q(.)$ defined on the interval $\Pq{t_0-\eps,t_0+\eps}$ is a minimizer between $q(t_0-\eps)$ and $q(t_0+\eps)$. For details, see e.g. \cite[Cor 17.1]{agra-book}. We apply this result to $q^1(t)$ in $t^1$ and find a corresponding $\eps^1$. Hence we have the minimizer $q^1(t)$ over the interval $\Pq{t^1-\eps^1,t^1+\eps^1}$. Notice that this minimizer is $\con^2$ but not $\con^3$ in $t^1$.

We now prove that for a small $\beps<\eps^1$ the trajectory $q^1(t)$ is not only a minimizer, but the unique \nuovo{normal} minimizer steering $Q_0=q^1(t^1-\beps)$ to $Q_1=q(t^1+\beps)$. We prove it by contradiction. Assume that there exists another minimizer $q^2$ steering $Q_0$ to $Q_1$. In the following all quantities related to this minimizer are denoted with superscript $2$. As a consequence of the existence of $q^2$, we have another minimizer $q^3(.)$ steering $q^1(t^1-\eps^1)$ to $q^1(t^1+\eps^1)$, given by the concatenation of $q^1$ on $\Pq{t^1-\eps^1,t^1-\beps}$, then $q^2$ on $\Pq{t^1-\beps,t^1+\beps}$, then again $q^1$ on $\Pq{t^1+\beps,t^1+\eps^1}$. See Figure \ref{fig:unicomin}. Since $q^3$ is a minimizer, then it is a solution of PMP. As a consequence, its tangent covector is continuous. For this reason, we have  $\lam^1(t^1-\beps)=\lam^2(t^1-\beps)$. Since this covector satisfies $h_1<0$, then trajectory $q^3$ satisfies dynamics given by \r{eq:dynbrutta} on a neighborhood of $t^1-\beps$, hence $q^1$ and  $q^2$ coincide on this neighborhood due to uniqueness of solution for \r{eq:dynbrutta}. We can prove in the same way that $q^1$ and $q^2$ coincide on the whole interval $[t^1-\beps,t^1)$. Similarly, we have $\lam^1(t^1+\beps)=\lam^2(t^1+\beps)$, hence $q^1$ and $q^2$ coincide in the whole interval $(t^1,t^1+\beps]$ due to uniqueness of solution for \r{eq:dynbella}. Finally, they coincide also in $t^1$ due to continuity. Hence $q^2=q^1$ on the interval $\Pq{t^1-\beps,t^1+\beps}$. Contradiction.
\immagine[9]{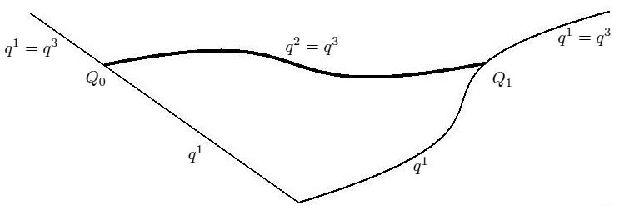}{Construction of trajectory $q^3$.}

We now define the trajectory $\q$ in $SE(2)$, using $q^1$ defined on the interval $\Pq{t^1-\beps,t^1+\beps}$. We first perform a left multiplication of $q^1$ in order to have $q^1(t^1)=\Id$, then a time shift $\Pq{t^1-\beps,t^1+\beps}\mapsto\Pq{0,2\beps}$. The resulting trajectory is  $\q(t):=\Pt{q^1(t^1)}^{-1} q^1(t+t^1-\beps)$. We recall some properties of this trajectory that we will use in the following:
\bi
\i $\q(\beps)=\Id$.
\i $\q$ is the unique minimizer steering $\q(0)$ to $\q(2\beps)$.
\i $\q$ satisfies dynamics \r{eq:dynbrutta} on $\Pq{0,\beps}$ and \r{eq:dynbella} on $\Pq{\beps,2 \beps}$. 
\ei

\subsection{Projective Reeds-Shepp car with lenght-curvature cost}
\label{s:control-proj}

The Reeds-Shepp car is a car that can move forwards, backwards and rotate on itself. The set of configurations can be identified with a quotient of the group of motions of the plane $SE(2)/\simeq$ where $(x,y,\th)\simeq(x,y, \th+\pi)$. For a better comprehension we use the same notation used for $SE(2)$, omitting the identification. We also omit checks of good definitions of dynamics and cost given below.

The dynamics of the car is given by the following control system:
\bqn
&&\Pt{\ba{c}
\dot x\\
\dot y\\
\dot \th
\ea
}=u_1\Pt{\ba{c}
\cos(\th)\\
\sin(\th)\\
0
\ea}+u_2\Pt{\ba{c}
0\\
0\\
1
\ea}
\eqnl{eq:dyncontrollo-proj}
where $u_1,u_2\in L^1([0,T],\R)$.

Fix a starting point $q_0=\Pt{x_0,y_0,\th_0}$ and an ending point $q_1=\Pt{x_1,y_1,\th_1}$. We want to minimize the cost
\bqn
\Cs\Pq{q(.)}=\int_0^T \sqrt{u_1^2+ u_2^2}
\eqnl{eq:costocontrollo-proj} over all trajectories of \r{eq:dyncontrollo-proj} steering $q_0$ to $q_1$. Here the end time $T$ is fixed.

Also in this case, due to invariance by rototranslation of both the dynamics and the cost, we can study only problems starting from $\Id=\Pt{0,0,0}=\Pt{0,0,\pi}$ and we will do so all along the paper.

Controllability is a direct consequence of Rashevsky-Chow theorem (see e.g. \cite{agra-book}) for this problem, since the distribution $\mathrm{span}\Pg{\Pt{\cos(\th),\sin(\th),0},\Pt{0,0,1}}$ is bracket-generating.

\subsubsection{Computation of extremals}
\label{s:extrem-proj}

In this section we compute minimizers for the optimal control problem \r{eq:dyncontrollo-proj}-\r{eq:costocontrollo-proj}. We follow procedure presented in Sections \ref{s:ex-min-SE2}-\ref{s:PMP-SE2}.

First transform it in a minimal time problem where dynamics is given again by \r{eq:dyncontrollo-proj} and controls are bounded by $u_1^2+u_2^2\leq 1$. Following Section \ref{s:ex-min-SE2}, we prove that the problem admits a minimum for all pair of starting and ending points and we restrict ourselves to $L^\infty$ optimal controls. We then apply PMP, using its expression for minimal time problem. Since dynamics \r{eq:dyncontrollo-proj} on $SE(2)/\simeq$ coincides locally with dynamics \r{eq:dyncontrollo} on $SE(2)$, we have the same control-depending Hamiltonian 
\bqn
H(q,\lam,u)=\Pa{\lam,\dot q}=u_1 h_1 + u_2 h_2
\eqnl{eq:Ham-proj}
where $h_1=\lam_x\cos(\th)+\lam_y \sin(\th),\ h_2=\lam_\th$, and $\lam_x,\lam_y,\lam_\th$ are the components of the covector $\lam$ in the dual basis with respect to coordinates $(x,y,\th)$. 

We can neglect abnormal extremals, since in this case they are trajectories  reduced to a point.
 
We fix $H=1$ and observe that we don't have condition $u_1\geq 0$ in this case. Hence solutions of the PMP are given by the choice $u_1=h_1,\, u_2=h_2$, that correspond to pendulum oscillations presented in Section \ref{s:PMP-SE2}. The corresponding dynamical system is the solution of
\bqn
\begin{cases}
\dot{x}=h_1 \cos(\th)\\
\dot{y}=h_1 \sin(\th)\\
\dot{\lam_x}=\dot{\lam_y}=0\\
\dot{\th}=h_2\\
\dot{\lam_\th}=h_1\Pt{-\lam_x\sin(\th)+\lam_y\cos(\th)}.
\end{cases}
\eqnl{eq:dynbella-proj}

The explicit solution of this problem is given in \cite{yuri-solo} in the case of $SE(2)$. For our treatment, it is sufficient to observe some properties of extremals. First of all, they are completely determined by the initial covector $\lam$, due to uniqueness of solution of \r{eq:dynbella-proj}. Moreover, the solution is analytic. As a consequence, we have only one of these possibilities:
\bi
\i either $h_1\equiv 0$, and the corresponding extremals are $q(t)=(0,0,\th^a(t))$,
\i or $h_1$ has only a finite number of times $t_1,\ldots,t_n$ in which it is vanishing, hence the corresponding trajectory $q(.)$ has only a finite number of points in which both $\dot{x}$ and $\dot{y}$ are vanishing.
\ei

Notice that trajectories of the second kind can be \virg{well projected} to the plane, i.e. it holds
\bl
\label{l:qb-bello}
Let $q(t)=(x(t),y(t),\th(t))$ be an extremal for the optimal control problem \r{eq:dyncontrollo-proj}-\r{eq:costocontrollo-proj} for which $h_1$ is vanishing only for a finite number of times $t_1,\ldots,t_n$. Let $p(t)=\Pi(q(t))$ be the projection of $q$ on the plane via $\Pi:(x,y,\th)\mapsto(x,y)$.  Then for each time $t\in[0,T]$ we have either $\dot{p}(t)\approx (\cos(\th(t)),\sin(\th(t)))$ or $\dot{p}(t_0)=0$ and $\lim_{\tau\rightarrow t^{-}} \frac{\dot{p}(\tau)}{\|\dot{p}(\tau)\|}\approx\lim_{\tau\rightarrow t^{+}} \frac{\dot{p}(\tau)}{\|\dot{p}(\tau)\|}\approx(\cos(\th(t)),\sin(\th(t)))$.
\el
\proof
First notice that $\dot{p}=(u_1 cos(\th),u_1 \sin(\th))$ since $q$ satisfies \r{eq:dyncontrollo-proj}. Hence it is clear that $\dot{p}(t)\approx (\cos(\th(t)),\sin(\th(t)))$ if $u_1(t)\neq 0$.

If instead $u_1(t)=0$, then there exists an interval $\Pt{t-\eps,t+\eps}$ on which $u_1(\tau)\neq 0$ for all $\tau\neq t$. Thus $\frac{\dot{p}(\tau)}{\|\dot{p}(\tau)\|}=\frac{(u_1(\tau) cos(\th(\tau)),u_1(\tau) \sin(\th(\tau)))}{|u_1(\tau)|}\approx(\cos(\th(\tau)),\sin(\th(\tau)))$. Passage to limit provides the result in $t$.
\qed

\brem An interesting property (see \cite{yuri-solo}) of this second family of extremals is that there are minimizers with one or two points in which $u_1=0$, but trajectory with three or more points in which $u_1=0$ are never minimizers. Thus minimizers for $J$ over the set $\dom_3$ may present one or two cusps, but not more than two.
\erem

We will use in the following these properties to prove existence of a minimizer of $J$ over all curves in $\dom_3$. Notice that Lemma \ref{l:qb-bello} doesn't hold for minimizers of the problem on $SE(2)$ defined in Section \ref{s:control-or}, since there are minimizers (like $\q$) such that their projection satisfies $\dot{\p}=0$ on an interval and $\dot{\p}\neq 0$ on another interval.

\section{Solution of problems and existence of minimizers}
\label{s:solution}

This section contains the main results of the paper. We first prove Propositions \ref{p:n-ex-1} and \ref{p:n-ex-2}, i.e. the non-existence of minimizers for the problem of minimization of $J$ respectively in $\dom_1$ and $\dom_2$. On the contrary, we prove Proposition \ref{p:ex-3}, that is the existence of minimizer for the problem of minimization of $J$ in $\dom_3$.

\subsection{Boundary conditions with orientation: non-existence of minimizers}
\label{ex:n-ex-1}
\label{ex:n-ex-2}

In this section we give a counterexample to the existence of minimizers of $J$ for boundary conditions with orientation. We prove it both in the case in which curves are chosen to be in $\dom_1$ and in $\dom_2$. The basic tools we use are the lift of a planar curve to $SE(2)$, see Remark \ref{r:lift}, and the trajectory $\b{q}(t)$ on $SE(2)$ defined in Section \ref{s:solstrana}, that is a solution of the optimal control problem \r{eq:dyncontrollo}-\r{eq:costocontrollo} studied in Section \ref{s:control-or}.

The basic idea is that we can lift the planar problem to the problem on $SE(2)$, then solve the problem on $SE(2)$ and finally project it again on the plane. But this last step doesn't work well, since in the case we present below the projection of the solution of the problem on $SE(2)$ doesn't satisfy boundary conditions with orientation fixed at the beginning.

Start considering the trajectory $\q(t)=(\b{x}(t),\b{y}(t),\b{\th}(t))$ on $SE(2)$ defined in Section \ref{s:solstrana} on the interval $\Pq{0,2 \beps}$. Define its projection $\p(t):=\Pi\Pt{\q(t)}$ on the plane $\R^2$ via the map $\Pi:(x,y,\th)\mapsto (x,y)$. As already stated, notice that $\dot{\p}=0$ on $\Pt{0,\beps}$ and $\dot{\p}\neq 0$ on $\Pt{\beps,2\beps}$. Then define a sequence of planar curves $p^n$ on the same interval, satisfying the following conditions:
\bi
\i Each of them satisfies the following boundary conditions with orientation: $$p^n(0)=\p(0),\, p^n(2\beps)=\p(2\beps),\, \dot p^n(0)\sim (\cos(\b{\th}(0)),\sin(\b{\th}(0))),\, \dot{p}^n(2\beps)\sim\dot{\b{p}}(2\beps).$$
\i The sequence converges to $\p$.
\i The cost $J\Pq{p^n}$ converges to $\Cs\Pq{\q}$.
\ei
Since now, notice that, if $p^n$ exists, then it is an example of the fact that each $p^n$ satisfies some boundary conditions with orientation but the limit trajectory $\p$ doesn't, since $\dot{\p}=0$ on $\Pt{0,\beps}$ and $\dot{\p}(\beps^+)\sim (1,0)$.

We define the curve $p^n$ with a geometric construction, see Figure \ref{fig:pn}. First define $p^n$ on $\Pq{\beps+\frac{\beps}{n},2 \beps}$ coinciding with $\p$. Then define the point $C:=\p\Pt{\beps+\frac{\beps}{n}}$ and draw the line $r$ that is the tangent to $\p$ or $p^n$ at $C$. Then draw the line $s$ passing through the origin $O=(0,0)$ and $(\cos(\b{\th}(0)),\sin(\b{\th}(0)))$. Since $\bth(\beps+t)=\bth(\beps)+t+o(t)=t+o(t)$, then $\bth(\beps+\frac{\beps}{n})>0$, while $\b{\th}(0)<0$, hence $r$ and $s$ are not parallel, thus they have an intersection point $B$.
\immagine[9]{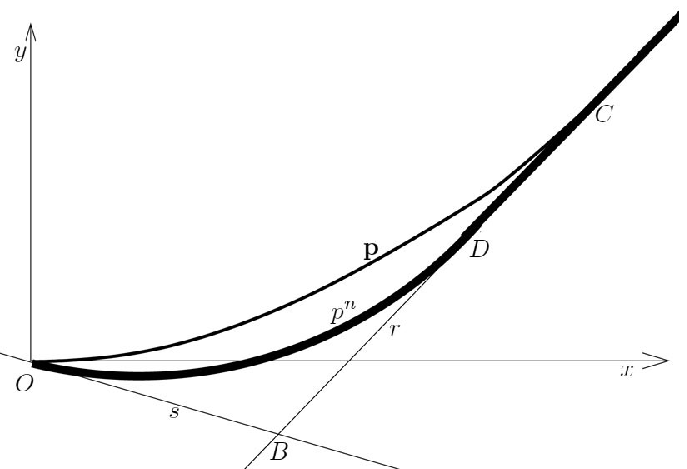}{Construction of the trajectory $p^n$ (case $\l{OB}\leq\l{BC}$).}

Then we have two cases:
\bi
\i If $\l{OB}\leq\l{BC}$, fix point $D$ on $BC$ such that $\l{OB}=\l{BD}$ and define the arc $\wideparen{OD}$ that is tangent to $OB$ in $O$ and to $BC$ in $D$. In this case, define $p^n$ on $\Pt{0,\beps+\frac{\beps}{n}}$ as the concatenation of the arc $\wideparen{OD}$ and the segment $DC$.

\i If instead $\l{OB}\geq\l{BC}$, fix $D$ on $OB$ satisfying $\l{BD}=\l{BC}$ and make the construction of the arc $\wideparen{DC}$. In this case $p^n$ on $\Pt{0,\beps+\frac{\beps}{n}}$ is the concatenation of the segment $OD$ and the arc $\wideparen{DC}$.
\ei

Notice that all $p^n$ satisfy boundary conditions with orientations and that the sequence converges to $\p$. Moreover, $J\Pq{p^n}$ restricted to the interval $\Pq{\beps+\frac{\beps}{n},2 \beps}$ coincides with $J\Pq{\p}$ on the same interval, that in turns coincide with $\Cs\Pq{\q}$ on the same interval, since $\q$ is the lift of $\p$ (see Remark \ref{r:lift}).

Concerning the interval $\Pq{0,\eps+\frac{\eps}{n}}$, we have $C\rightarrow B\rightarrow O$ for $n\rightarrow\infty$, hence $J\Pq{DC}$ or $J\Pq{OD}$ tends to 0. Instead the cost of the arc $\wideparen{OD}$ or $\wideparen{DC}$ tends to $-\bth(0)$. Indeed, assume that $\l{OB}\leq\l{BC}$ and compute $J\Pq{\wideparen{OD}}$ with an arclength parametrization. Recall that in this case $\int_a^b K_\ga=\al_\ga(b)-\al_\ga(a)$ where $\al_\ga(t)$ is the angle of the tangent vector $\dot\ga$, thus  $$\al_{\wideparen{OD}}\Pt{D}-\al_{\wideparen{OD}}(O)=\int_0^{\l{\wideparen{OD}}} K_\ga\,ds\leq J\Pq{\wideparen{OD}}\leq\int_0^{\l{\wideparen{OD}}} (1+K_\ga)\,ds=\l{\wideparen{OD}}+\al_{\wideparen{OD}}\Pt{D}-\al_{\wideparen{OD}}(O).$$
The result follows recalling that $\al_{\wideparen{OD}}\Pt{D}\rightarrow \bth\Pt{\beps+\frac{\beps}{n}}\rightarrow \bth\Pt{\beps}=0$. The case $\l{OB}\geq\l{BC}$ can be treated similarly.

We have thus defined a sequence of curves $p^n\in\dom_1\subset\dom_2$ minimizing the cost $J$ but such that the limit curve $\p$ does not satisfy boundary conditions with orientation, hence it is not in $\dom_2$. We prove that it implies the non-existence of a minimizer for these boundary conditions. By contradiction, assume that a minimizer of $J$ exists in $\dom_2$. Thus its lift $\bar{q}$ to $SE(2)$ is a minimizer for $\Cs$ between $\q(0)$ and $\q(2\eps)$. Since $\q$ is the unique \nuovo{normal minimizer between the two points, then $\bar{q}$ is abnormal. Since $\Pt{{\bf x}(0),{\bf y}(0)}\neq \Pt{{\bf x}(2\eps),{\bf y}(2\eps)}$ and ${\bf \th}(0)\neq {\bf \th}(2\eps)$, then $\bar{q}$ is neither a straight line nor a pure rotation, hence it is a concatenation of straight lines and rotations. Its projection is thus a curve with angles, i.e. it is not in $\dom_2$. Contradiction.}

\subsection{Projective boundary conditions with orientation: existence of minimizers}
\label{s:sol-exist}

In this section we prove existence of a minimizing curve in $\dom_3$ for all choices of projective boundary conditions, i.e. we prove Proposition \ref{p:ex-3}. The basic idea is that also in this case we can lift the problem of planar curves to the problem on $SE(2)/\simeq$ defined above, solve it and then project the solution to the plane. But in this case the whole procedure works well, since the projection of the solution of the problem on $SE(2)/\simeq$ is always the solution of the planar problem. In particular, it satisfies projective boundary conditions.

Start fixing projective boundary conditions, i.e. fix a starting point $(x_0,y_0)$ with direction $v_0$ and an ending point $(x_1,y_1)$ with direction $v_1$. Assume that $(x_0,y_0)\neq (x_1,y_1)$ and $v_0,\, v_1$ are nonvanishing vectors. Recall that we want to find a curve $\ga\in\dom_2$ such that $\ga(0)=(x_0,y_0)$, $\dot\ga(0)\approx v_0$, $\ga(T)=(x_1,y_1)$, $\dot\ga(T)\approx v_1$ and that is a minimizer of $J$. 

Consider the optimal control defined on $SE(2)/\simeq$ presented in Section \ref{s:control-proj} with the following starting and ending point: $q_0=(x_0,y_0,\th_0)$ and $q_1=(x_1,y_1,\th_1)$ where each $\th_i$ is the angle formed by the vector $v_i$ with respect to the $x$-axis. Then solve the problem and call $\q(.)$ the minimizing trajectory (that is not necessary unique). The basic remark is that $\q$ is of the second kind (see Section \ref{s:extrem-proj}), since $(x_0,y_0)\neq (x_1,y_1)$. As proved in Lemma \ref{l:qb-bello}, in this case $\dot{\p}\approx (\cos(\bth),\sin(\bth))$ except for a discrete set of points $t_1,\ldots, t_n$ on which we have the weaker property $\lim_{\tau\rightarrow t_i}\frac{\dot{\p}}{\|\dot{\p}\|}\approx (\cos(\bth),\sin(\bth))$. If we have at the starting point $\dot{\p}(0)\neq 0$, then $\p$ satisfies projective boundary conditions at the beginning. Otherwise reparametrize $\p$ by arclength in an interval $\Pq{0,\eps}$, that is possible since 0 is the unique point in the interval in which $\dot{\p}=0$. As a consequence, now $\p$ satisfies boundary conditions at the beginning, since we have $\dot{\p}(0)=\lim_{\tau\rightarrow 0}\frac{\dot{\p}}{\|\dot{\p}\|}\approx (\cos(\bth(0)),\sin(\bth(0)))$. The same result can be proved for the ending point. Hence $\p$ satisfies projective boundary conditions.

We now prove that $\p$ is a minimizer of $J$, by contradiction. Assume that there exists $\bar{p}$ satisfying the same projective boundary conditions and such that $J\Pq{\bar{p}}<J\Pq{\p}$. Thus its lift $\bar{q}$ steers $q_0$ to $q_1$ and satisfies $\Cs\Pq{\bar{q}}<\Cs\Pq{\q}$, hence $\q$ is not a minimizer. Contradiction.

\end{document}